\documentclass[a4paper]{jpconf}
\usepackage{graphics}
\usepackage{epsfig}
\usepackage{amssymb,amsmath}
\usepackage{amsfonts}
\usepackage{color}
\usepackage{pstricks}
\usepackage{epstopdf}

\catcode`@=11 \@addtoreset{equation}{section}
\renewcommand\theequation{\thesection.\@arabic\c@equation}
\catcode`@=12

\makeatother

\makeatother

\begin{document}
\title{A Two-Layer Mathematical Modelling of Drug Delivery to Biological Tissues}

\author{Koyel Chakravarty and D C Dalal}

\address{Department of Mathematics, Indian Institute of Technology Guwahati, Guwahati, Assam, India}

\ead{koyel@iitg.ernet.in}

\begin{abstract}
Local drug delivery has received much recognition in recent years, yet it is still unpredictable how drug efficacy depends on physicochemical properties and delivery kinetics. The purpose of the current study is to provide a useful mathematical model for drug release from a drug delivery device and consecutive drug transport in biological tissue, thereby aiding the development of new therapeutic drug by a systemic approach. In order to study the complete process, a two-layer spatio-temporal model depicting drug transport between the coupled media is presented. Drug release is described by considering solubilisation dynamics of drug particle, diffusion of the solubilised drug through porous matrix and also some other processes like reversible dissociation / recrystallization, drug particle-receptor binding and internalization phenomena. The model has led to a system of partial differential equations describing the important properties of drug kinetics. This model contributes towards the perception of the roles played by diffusion, mass-transfer, particle binding and internalization parameters.
\end{abstract}

\section{Introduction}\label{secint}
Studies related to drug delivery draw much attention to present day researchers for their theoretical and clinical investigations.  The controlled drug delivery is a process by means of which a drug is delivered at a pre-determined rate, locally or systemically, for a stipulated period of time. Controlled drug delivery systems may incorporate the maintenance of drug levels within a desired therapeutic range, the need of less number of administrations, optimal use of the required drug and yet possibility of enhanced patient compliance. The quintessential drug delivery should be not only inert, biocompatible, but also at the same time provide patient compliance, capable of attaining high drug loading. However, it should also have preventive measures for the accidental release of drug and at the same time it should be simple to administer and to remove from the body. Controlled release drug delivery involves drug-encapsulating devices from which drug or therapeutic compounds may be released at controlled rates for prolonged span of time. \\

While many old and new therapeutics are well tolerated, numerous compounds are in requirement of localized advanced drug delivery technologies to reduce toxicity level, enhance therapeutic efficacy and potentially recast bio-distribution. Local drug delivery is the manifestation of drug delivery in which drugs are delivered at a specific site inside the body to a particular diseased organ or tissue. Though the drug delivery, in principle, may be monitored, but the most important hazard is that the design of drug delivery is unclear, which must be used to attain the level of control required for a specific purpose. This is because there exists complex interaction between biology, polymer chemistry and pharmacology \cite{1}.\\

Mathematical modelling of drug delivery and predictability of drug release is a steadily growing field with respect to its importance in academic and industrial areas due to its astronomic future potential. In the light of drug dose to be incorporated in desired drug administration and targeted drug release profile, mathematical prognosis will allow for good estimates of the  required composition along with other requirements of the respective drug dosage forms. An extremely challenging aspect is to combine mathematical theories with models quantifying the release and transport of drug in living tissues and cells. Various works are done in the past on drug delivery devices regarding its therapeutic efficiency, optimal design with the aid of either experimental methods or numerical / modelling simulations and sometimes both procedures are used \cite{2,3,4,5}. The investigation in which drug association / dissociation aspect is taken into account with regard to transdermal drug delivery \cite{6} while in another study, drug release process from a microparticle is considered based on solubilisation dynamics of the drug \cite{7}. Very recently, the updated mathematical model considering both the above mentioned aspects is framed and analysed successfully \cite{8}.\\

 The study of concern is presented by a phenomenological mathematical model of drug release from a local drug delivery device and its subsequent transport to the biological tissue. The model is comprised of two-phase drug release where the drug undergoes solubilisation, recrystallisation and internalization through a porous membrane. An important aspect of the aforementioned model is appropriate judgement of the model parameters of significance such as diffusion coefficient, solid-liquid transfer rate, mass-transfer, drug association / dissociation rate constants, membrane permeability and internalization rate constant. The numerical simulation provides reliable information on different properties of drug release kinetics. \\


\section{Problem of drug release}
To model the local drug delivery device, a two-phase system is considered which is made of: (a) a polymeric matrix that operates as a reservoir where the drug is loaded initially, and (b) the biological tissue where the drug is being transported as target region. The first phase i.e the drug delivery device is framed as a planar slab, encompassed on one side with an impermeable backing and the other side of the device is in contact with layer (b). A rate-controlling membrane protecting the polymeric matrix is present at the interface of the coupled layers.\\ 

\subsection{\bf{Mechanism of drug release }}
At the beginning, the drug occurs completely in a solid phase embraced within the polymeric matrix (e.g. in crystalline form) ($C_0^{*}$) at its maximum concentration. Being in bound state, it cannot be transferred to the tissues directly. Water enters into the polymeric matrix and wets the drug encapsulated inside it, permitting solubilisation of the loaded drug crystals into free state ($C_0$) which diffuses out of the matrix into the tissue. The rate of transfer of drug from solid state to free state depends not only on solubilisation phenomenon but is also proportional to the difference between $C_0^{*}$ and $C_0$. Again, a fraction of solid drug ($\beta_0C_0^{*}$) is transformed to its free state which is competent to diffuse. Conversely, through a recrystallisation process, another fraction of free drug ($\delta_0C_0$) is transferred back to its bound state. Simultaneously, a part of free drug ($C_1$) diffuses into the tissue. In the similar way, in tissue, a portion of free drug ($k_aC_1$) is metabolised into bound phase ($C_1^{*}$), which also unbinds ($k_dC_1^{*}$) to form free drug. Now, the bound drug is engulfed (internalized) ($k_iC_1^{*}$) by the cell in the tissue through the process of endocytosis. Endocytosis is an active energy-using transport phenomenon in which molecules (proteins, drugs etc.) are transported into the cell. Thus, bound drug gets transformed into internalized drug particles ($C_i$). These internalized drug particles, after a span of time, gets degraded by the lysosomes and the drug remnants after degradation ($k_{id}C_i$) is expelled out of the cell into the extracellular fluid. The complete drug transport process is schematically demonstrated in Fig. 1.  
\begin{figure}[ht!]
\centering
\includegraphics[width=14cm,height=3cm]{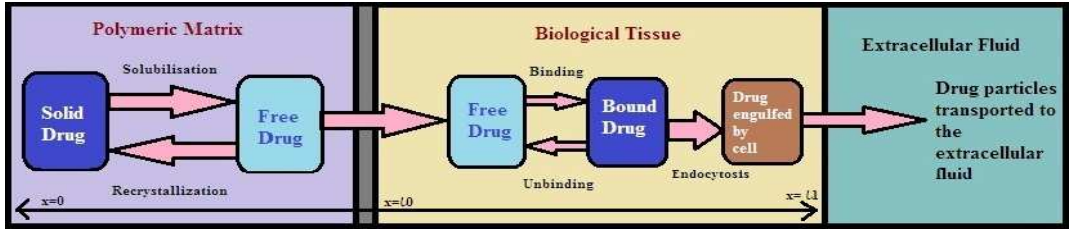}
\caption{\label{figure 1}  Schematic diagram of drug transport}
\end{figure}

\section{Formulation of the model}
Generally, mass transport prevails in the direction normal to the tissue which may be the reason behind the modelling restriction confined to one-dimensional case. In the present study, x-axis is considered to be normal to the layer and aligned  with the positive direction outwards.\\
\subsection{\bf{Drug dynamics in the polymeric matrix}}
The governing equations describing the dynamics of drug release in the polymeric matrix phase are
\begin{eqnarray}
\frac{\partial C_0^*}{\partial t}=-\alpha_0(\Phi_0 C_0^*-C_0)-k_m(C_{lim}-C_0)-\beta_0C_0^*+\delta_0C_0,      \qquad x\in (0,l_0)\\
\frac{\partial C_0}{\partial t}= D_0\frac{\partial^2 C_0}{\partial x^2}+\alpha_0(\Phi_0C_0^*-C_0)+k_m(C_{lim}-C_0)+\beta_{0}C_{0}^*
-\delta_0C_0, \qquad x\in (0,l_0)
\end{eqnarray}
where,
$\Phi_0 (=\frac{k\varepsilon_0}{1-\varepsilon_0})$
is the ratio of accessible void volume to solid volume,
$C_0^*$ denotes available molar concentration of solid drug,
$C_0$ is the available molar concentration of free drug,
$k$ stands for partition coefficient,
$\varepsilon_0$ denotes porosity,
$l_0$ is the length of the polymeric matrix,
$k_m$ is mass transfer coefficient,
$C_{lim}$ stands for drug solubilisation limit,
$\beta_0$ is the dissociation rate constant,
$\delta_0$ is the association rate constant,
$\alpha_0$ denotes solid-liquid rate parameter.
$D_0$ is the diffusion coefficient of free drug in the matrix.\\
\subsection{\bf{Drug dynamics in the biological tissue}}
The corresponding equations governing the dynamics of drug in the tissue are
\begin{eqnarray}
\frac{\partial C_1^*}{\partial t}=k_aC_1-k_dC_1^*-k_iC_1^*,      \qquad x\in (l_0,l_1)\\
\frac{\partial C_1}{\partial t}=D_1\frac{\partial^2 C_1}{\partial x^2}-k_aC_1+k_dC_1^*,      \qquad x\in (l_0,l_1)\\
\frac{\partial C_i}{\partial t}=k_iC_1^*-k_{id}C_i,                 \qquad x\in (l_0,l_1)
\end{eqnarray}
where,
$C_1$ is the available molar concentration of free drug in the tissue,
$C_1^*$ is the available molar concentration of solid drug in the tissue,
$C_i$ denotes the molar concentration of internalized drug particles,
$l_1-l_0$ is the length of the tissue,
$k_a$ depicts the binding rate coefficient,
$k_d$ is the dissociation rate coefficient,
$k_i$ stands for internalization rate coefficient,
$k_{id}$ denotes degradation rate constant in the lysosome, and
$D_1$ is the diffusion coefficient of free drug in the biological tissue.\\

\section{Initial, Interface and Boundary Conditions}
\begin{flushleft}
The initial conditions are as follows:\\
 $C_0^*(x,0)=M$, $C_0(x,0)=0$, $C_1^*(x,0)=0$, $C_1(x,0)=0$, $C_i(x,0)=0$.\\
  A flux continuity must be assigned at the interface, i.e at $x=l_0$, $-D_0\frac{\partial C_0}{\partial x}=-D_1\frac{\partial C_1}{\partial x}$.\\
No mass flux can pass to exterior environment due to the presence of impermeable backing and hence no flux condition arises.
At $x=0$, $D_0\frac{\partial C_0}{\partial x}=0$.\\
Lastly, at $x=l_1$, $C_1(l_1,t)$ is considered to be finite.\\
\end{flushleft}
\section{Model Solutions}

For the purpose of reducing the number of model parameters, the entire above-mentioned equations together with all the conditions are made dimensionless. The transformed equations subject to the conditions in dimensionless approach are solved analytically by separation of variables procedure. Thus, the solutions should be read as\\ 
\begin{equation}
C_0=E_1(e^{-m_1t}-e^{-m_2t})\cos ax\\
\end{equation}
\begin{equation}
C_0^{*}=\frac{E_1(\alpha_0+k_m+\delta_0)\cos ax}{(m_1-(\alpha_0\Phi_0+\beta_0))(m_2-(\alpha_0\Phi_0+\beta_0))}\Big[(m_1-(\alpha_0\Phi_0+\beta_0))e^{-m_2t}-(m_2-(\alpha_0\Phi_0+\beta_0))\\
\end{equation}
\begin{equation*}
e^{-m_1t}+(m_2-m_1)e^{-(\alpha_0\Phi_0+\beta_0)t}\Big]
\end{equation*}
\begin{equation*}
-\frac{k_mC_{lim}}{(\alpha_0\Phi_0+\beta_0)}\Big[1-e^{-(\alpha_0\Phi_0+\beta_0)t}\Big]+e^{-(\alpha_0\Phi_0+\beta_0)t}
\end{equation*}
\begin{equation}
C_1=E_2(e^{-n_1t}-e^{-n_2t})\cos bx\\
\end{equation}

\begin{equation}
C_1^{*}=\frac{E_2k_a \cos bx}{(k_d+k_i-n_1)(k_d+k_i-n_2)}\Big[(k_d+k_i-n_2)e^{-n_1t}-(k_d+k_i-n_1)e^{-n_2t}+(n_2-n_1)e^{-(k_d+k_i)t}\Big]
\end{equation}
\begin{equation}
C_i=\frac{E_2k_a k_i \cos bx}{(k_d+k_i-n_1)(k_d+k_i-n_2)}\Big[\frac{(k_d+k_i-n_2)}{(k_{id}-n_1)}(e^{-n_1t}-e^{-k_{id}t})-\frac{(k_d+k_i-n_1)}{(k_{id}-n_2)}(e^{-n_2t}-e^{-k_{id}t})+\\
\end{equation}
\begin{equation*}
\frac{(n_2-n_1)}{(k_{id}-k_d-k_i)}(e^{-(k_d+k_i)t}-e^{-k_{id}t})\Big]
\end{equation*}
where, $E_1$ and $E_2$ are arbitrary constants to be determined from the prevailing conditions,
\begin{eqnarray*}
m_1=\frac{A-\sqrt{A^2+4B}}{2},\qquad m_2=\frac{A+\sqrt{A^2+4B}}{2}\\
A=\alpha_0(\Phi_0+1)+k_m+\beta_0+\delta_0-\lambda\gamma, \qquad B=\lambda\gamma(\alpha_0\Phi_0+\beta_0),
\qquad \lambda=-a^2,\\
n_1=\frac{P-\sqrt{P^2-4Q}}{2},\qquad n_2=\frac{P+\sqrt{P^2-4Q}}{2}\\
P=k_d+k_i+k_a-\mu,\qquad Q=k_ik_a-\mu(k_d+k_i)
,\qquad \mu=-b^2.
\end{eqnarray*}
and, $a$ and $b$ are any positive real numbers.\\
All the parameters and the variables are expressed in dimensionless form whose expressions are not included here for the sake of brevity. 
\section{Numerical illustration and discussion}
Numerical illustration for the present drug release system is performed by taking various parameter values of the model in order to characterise the pharmacokinetic aspects.
The graphical representations of the time-variant concentration profiles of the drug in its different states for both the facets are well illustrated through the Figs. 2 - 6 in order to understand the drug release phenomenon.\\
The time-variant concentration profiles for solid (loaded) and free drug particles in the polymeric matrix phase for four different axial locations spread over the entire domain are demonstrated in Figs. 2 and 3.
\begin{figure}[h]
\begin{minipage}{3in}
\includegraphics[width=3in]{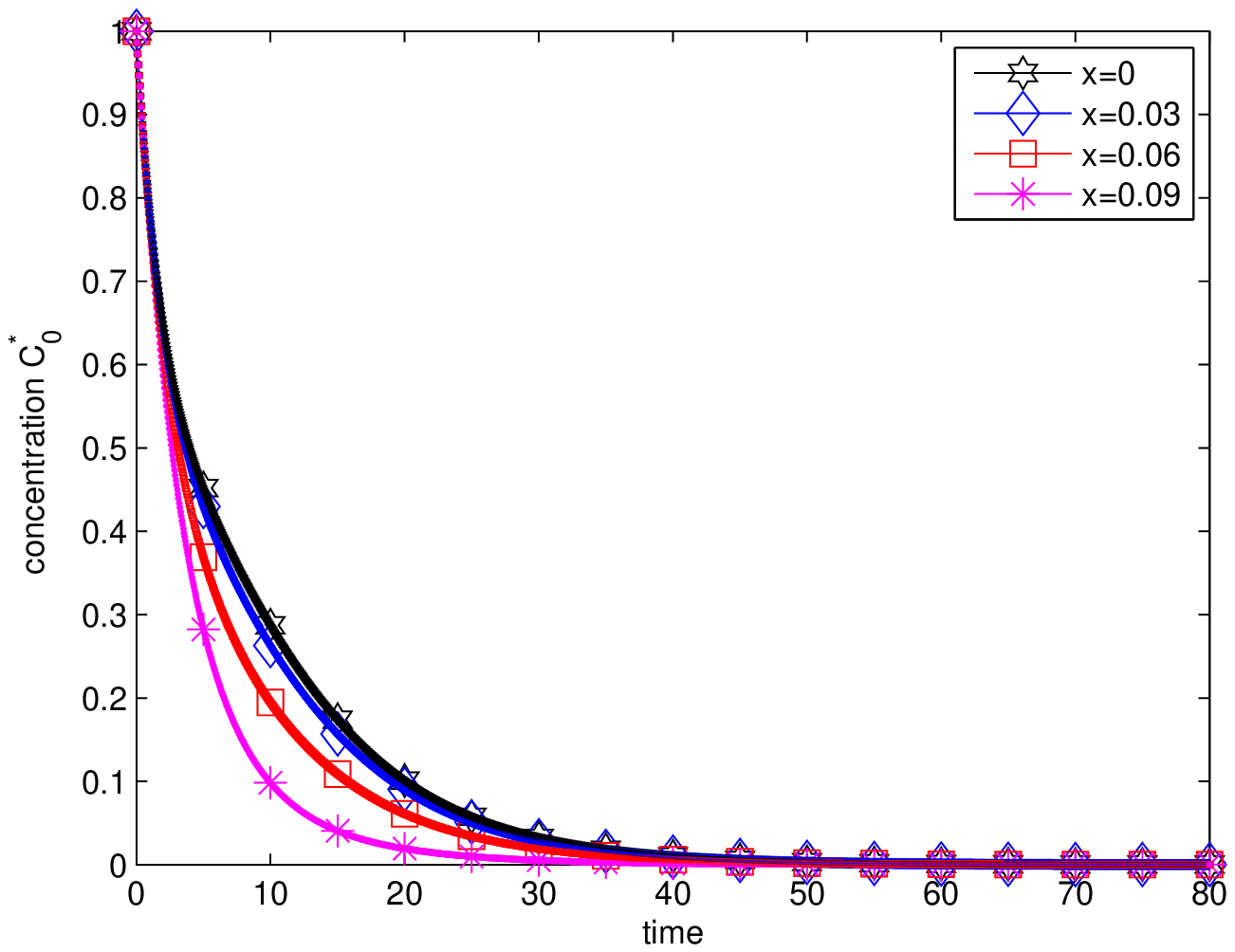}
\caption{\label{figure 2} Time variant concentration profile of $C_0^{*}$ at different locations.}
\end{minipage}
\hspace{.2in}
\begin{minipage}{3in}
\includegraphics[width=3in]{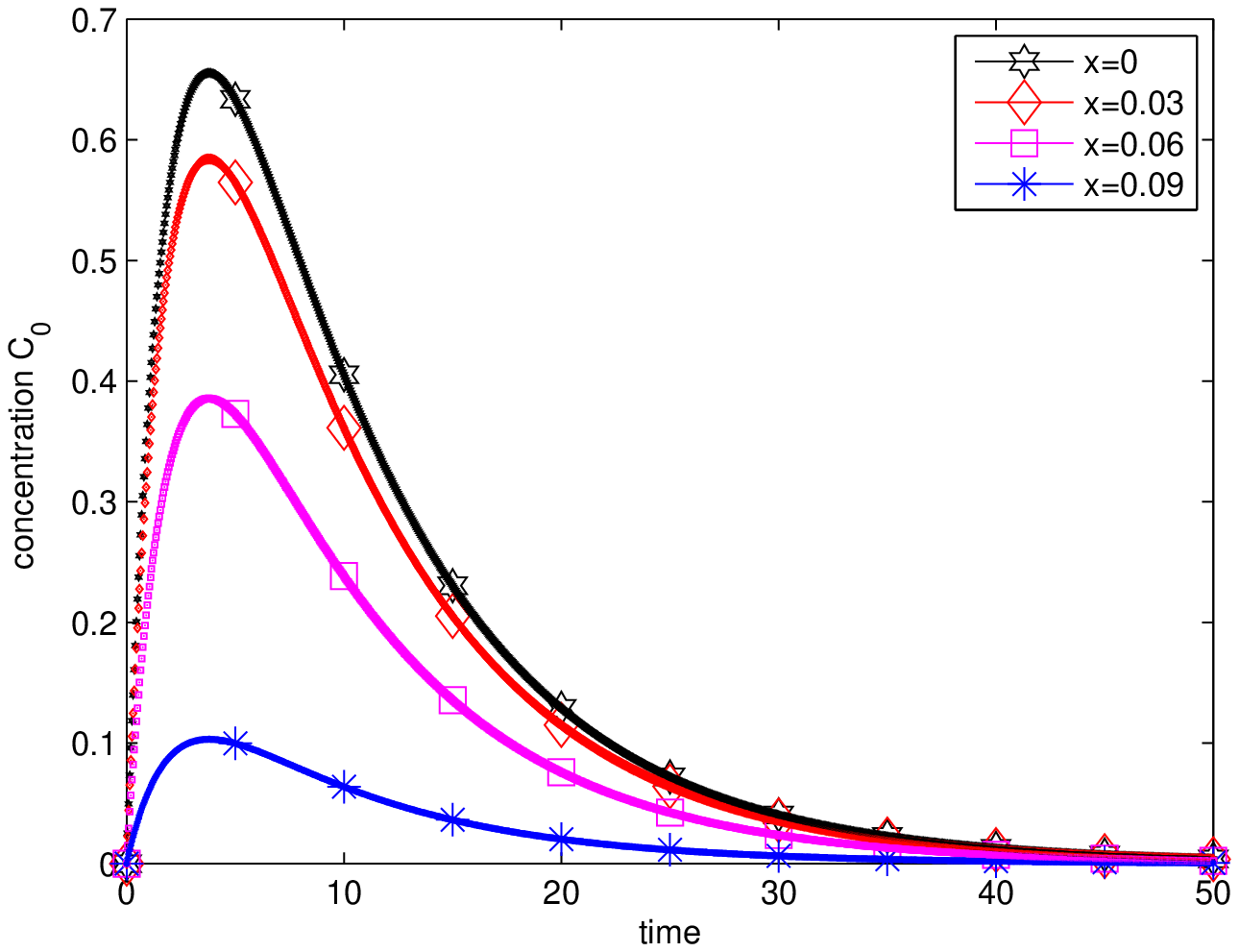}
\caption{\label{figure 3} Time variant concentration profile of $C_0$ at different locations.}
\end{minipage}
\end{figure}
 The rate of decrease of concentration of loaded drug ($C_0^{*}$) becomes higher and higher as one proceeds away from the commencing region of the polymeric matrix to the interface resulting in early disappearance. When the loaded drug gets exposed to water, a solid-liquid mass transfer is initiated causing drug release from the matrix which is on the process of transformation of free drug $C_0$. One may note on the other hand that $C_0$ grows and acquires certain peak in accordance to a specific instant of time followed by a gradual descend for rest of the times. \\
 \begin{figure}[h]
\begin{minipage}{2in}
\includegraphics[width=2.3in]{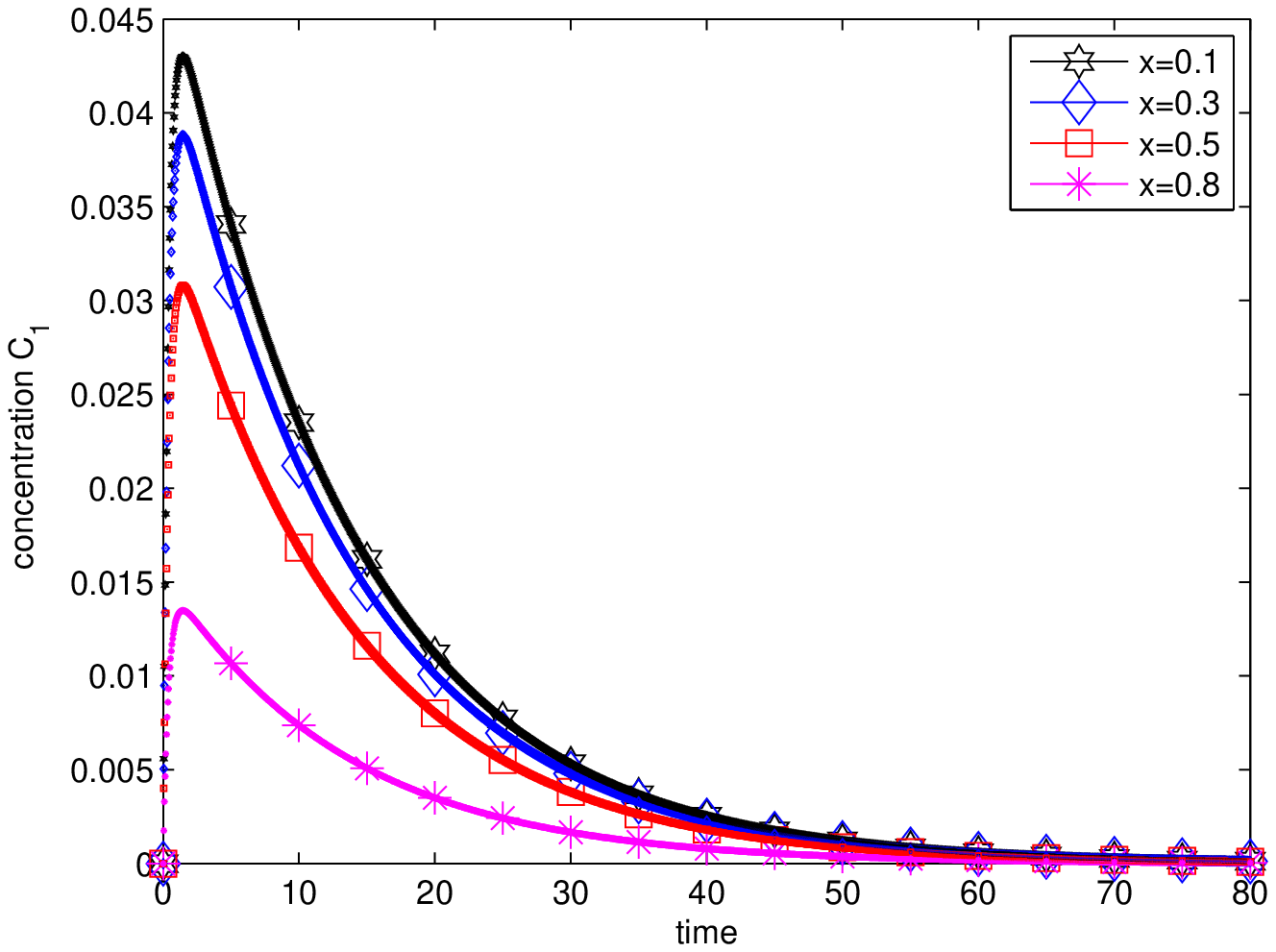}
\caption{\label{figure 4} Time variant concentration profile of $C_1$ at different locations.}
\end{minipage}
\hspace{.09in}
\begin{minipage}{2in}
\includegraphics[width=2.3in]{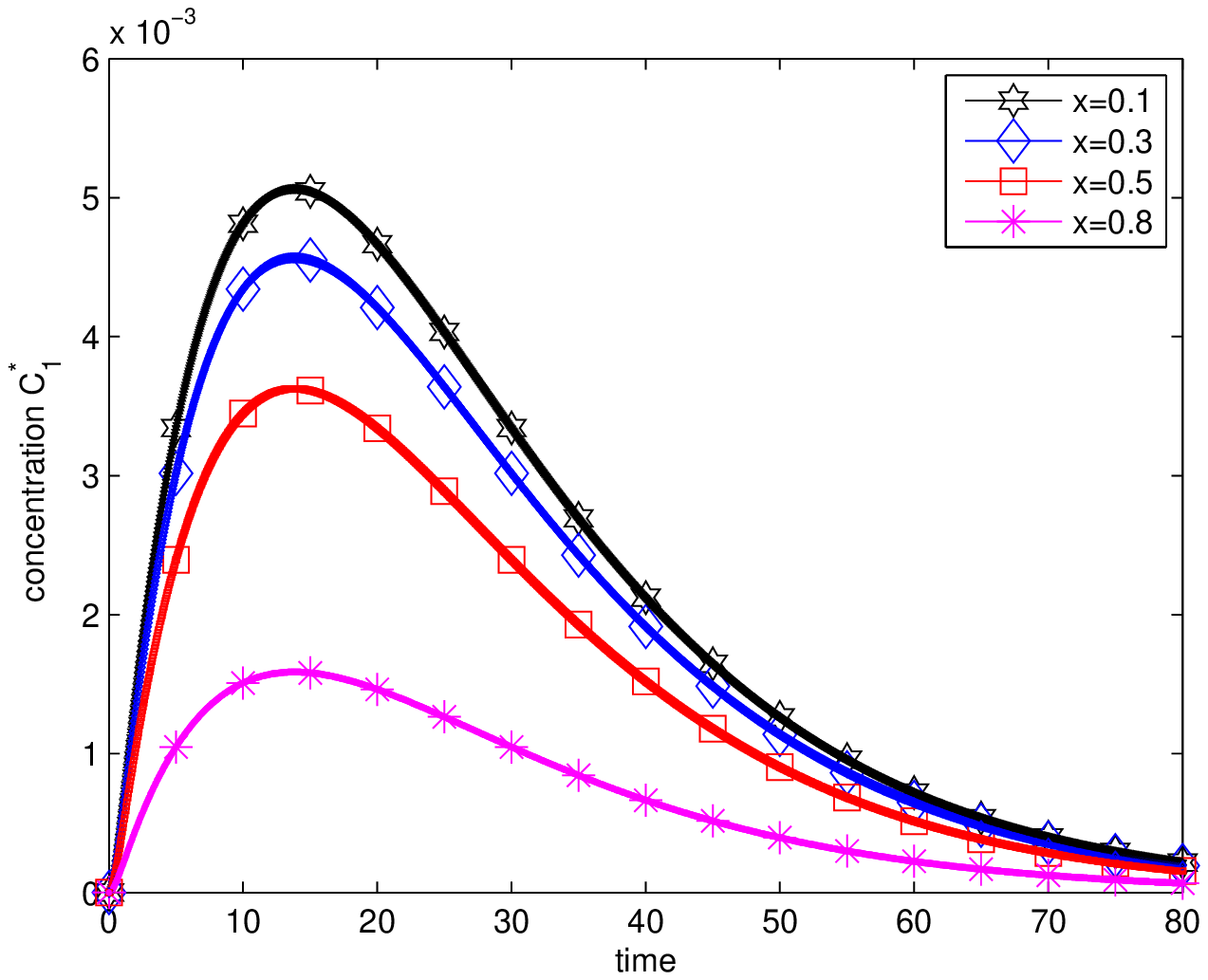}
\caption{\label{figure 5} Time variant concentration profile of $C_1^{*}$ at different locations.}
\end{minipage}
\hspace{.09in}
\begin{minipage}{2in}
\includegraphics[width=2.3in]{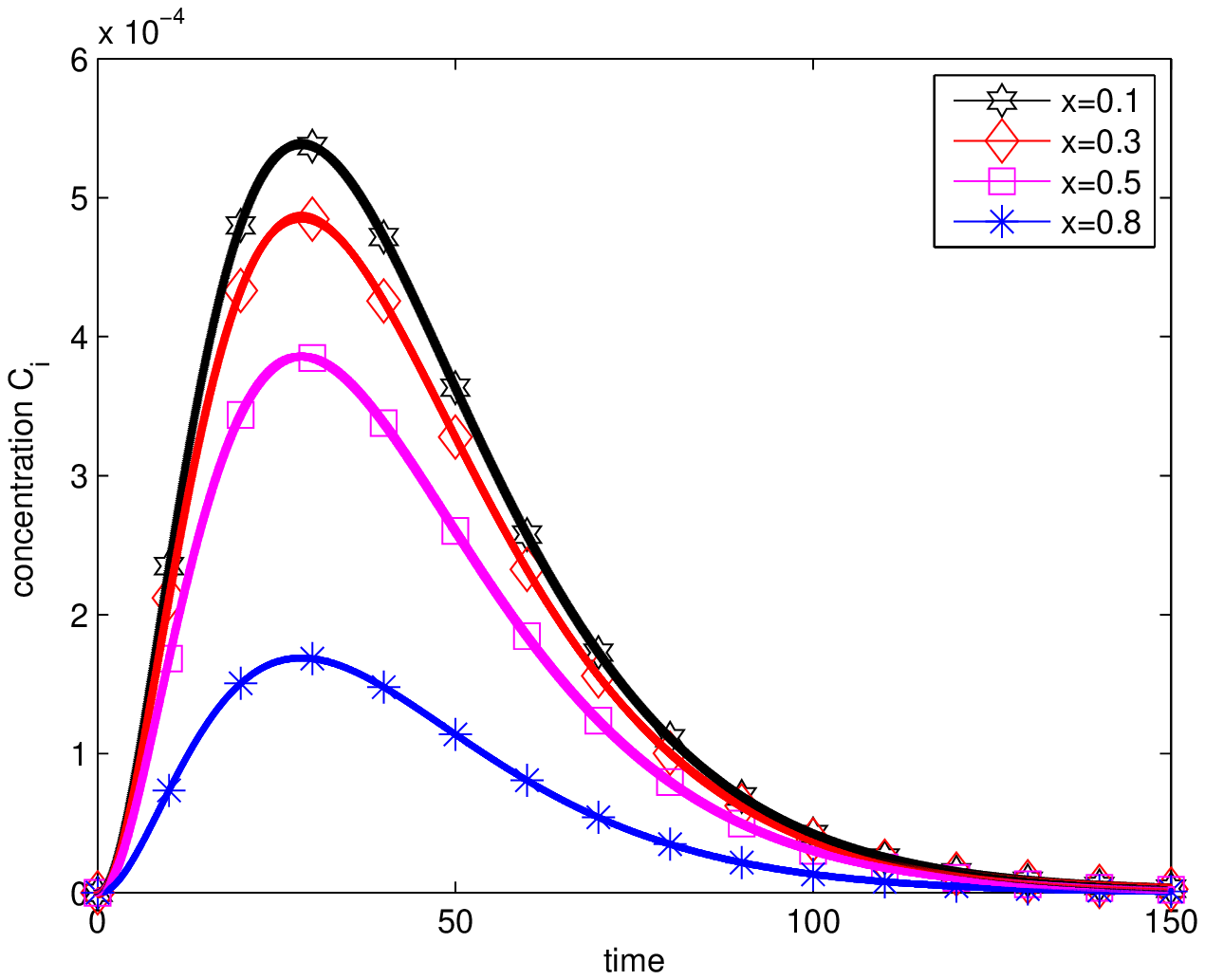}
\caption{\label{figure 6} Time variant concentration profile of $C_i$ at different locations.}
\end{minipage}
\end{figure}
Figs. 4 - 6 represent the time-variant concentration profiles for free, bound and internalized drug particles respectively in the biological tissue for different axial locations stretched over the entire domain. It is important to observe that $C_1$ heightens hastily compared to both $C_1^{*}$ and $C_i$ towards the inception. This observation, as anticipated, reflects in the realm of drug kinetics that free drug ($C_1$) gets transformed into bound drug ($C_1^{*}$) and subsequently the bound drug is metabolised into internalized drug ($C_i$) after a short passage of time. Ultimately, it is further noted that the internalized drug takes more time to get absorbed completely in the tissue than the characteristics of both bound and free drug particles. In addition to the present findings, one may append that the extended time span distinctly reveals that both loaded and free drug particles in polymeric matrix melt away in a comparatively small span of time than those in the tissue where they need to take time lengthened to get the drug absorbed completely.\\
One may also explore a variety of cases in order to exhibit the behaviour of the concentration profiles for both the phases under present consideration by varying all the parameter values of significance. The sensitivity of the model parameters imply the need of the components to be included in the model formation for future course of studies relevant to this domain.\\  
\section{Conclusions}
 The sensitivity of the model parameters poses challenges to the applicability of drug administration for treatment of patients at large through pharmacotherapy. One may highlight that as both loaded and free drug particles in polymeric matrix melt away in comparatively small span of time than those in the tissue, its influence will certainly persist for a long time before repeated medication occurs and hence care needs to be exercised for maintaining appropriate time-gap before redispensation  in order to avoid toxicity by the presence of excess drug.  \\


\end{document}